\documentclass{article}
\usepackage{amsmath,amsfonts,amssymb}
\usepackage[latin2]{inputenc}                                          
\def\softd{{\leavevmode\setbox1=\hbox{d}%
\hbox to 1.05\wd1{d\kern-0.4ex{\char039}\hss}}}
\def\softt{{\leavevmode\setbox1=\hbox{t}
\hbox to \wd1{t\kern-0.6ex{\char039}\hss}}}
\def\softl{l\kern-0.45ex\raise0.1ex\hbox{''}\kern-0.10ex}
\def\softL{L\kern-0.8ex\raise0.1ex\hbox{''}\kern0.1ex}

\newtheorem{theorem}{Theorem}[section]
\newtheorem{definition}[theorem]{Definition}

\newtheorem{example}[theorem]{Example}
\newtheorem{lemma}[theorem]{Lemma}

\newtheorem{remark}[theorem]{Remark}

\newcommand{\mC}{\mathbb C}

\newcommand{\lL}{\cal L}

\newcommand{\mZ}{\mathbb Z}

\newcommand{\mR}{\mathbb R}

\newcommand{\mN}{\mathbb N}
\newcommand{\be}{\begin{eqnarray}}
\newcommand{\ee}{\end{eqnarray}}
\newcommand{\bd}{\begin{definition}}
\newcommand{\ed}{\end{definition}}
\newcommand{\br}{\begin{remark}}
\newcommand{\er}{\end{remark}}

\newcommand{\bt}{\begin{tabular}}
\newcommand{\et}{\end{tabular}}

\newcommand{\kl}{{\fam2 L}}

\newcommand{\fn}{{\fam2 C}}

\newcommand{\bl}{\begin{lemma}}
\newcommand{\el}{\end{lemma}}

\newcommand{\bp}{\begin{picture}}
\newcommand{\ep}{\end{picture}}
\newcommand{\bi}{\begin{itemize}}
\newcommand{\ei}{\end{itemize}}
\newcommand{\bq}{\begin{quotation}}
\newcommand{\eq}{\end{quotation}}

\input epsf
\input xy
\xyoption{all}

\begin{document}
\baselineskip13pt

\title{Finite reflection groups and the 
Dunkl-Laplace differential-difference operators in conformal geometry
}

\author{P. Somberg}

\date {}

\maketitle

\abstract 
For a finite reflection subgroup $G\leq O(n+1,1,\mR)$ of the conformal group of
the sphere with standard conformal structure $(S^n,[g_0])$, we geometrically 
derive differential-difference Dunkl version
of the series of conformally invariant differential operators with symbols given by 
powers of Laplace operator. The construction can be regarded as a deformation of 
the Fefferman-Graham ambient metric construction of GJMS operators. 

\hspace{0.5cm}

{\bf Key words:} Conformal geometry, Finite reflection groups, Fefferman-Graham ambient 
metric construction, conformal Dunkl-Laplace operators.
 
{\bf MSC classification:} 51F15, 53A30, 20F55, 58E09.

\endabstract


\section{Introduction}

Invariant theory for finite reflection groups G is usually
considered in the framework $G\leq O(n)$, where $O(n)$ is 
the orthogonal group associated to the
Hilbert space $(\mR^n, <,>)$. The action of $G$ then carries over 
to function spaces on $\mR^n$ and leads to the concept of Dunkl 
differential-difference operators and special function theory
for finite reflection groups, see e.g. \cite{eting}, \cite{opdam}, 
\cite{ros} and
extensive references therein.

  The aim of the present article is to change the geometrical perspective and 
study similar questions (or at least, some of them) as in the classical 
well-known case. In particular, our point of view is to 
consider the conformal compactification $S^n$ of $\mR^n$ and finite 
reflection subgroups of the conformal group $O(n + 1, 1,\mR)$ defined
by finite number of reflecting subspheres diffeomorphic to spheres 
$S^{n-1}$ of codimension one. The collection of reflecting subspheres is 
an analog of the collection of reflecting hyperplanes in $\mR^n$.
Moreover, a reflection subgroup of the conformal group is induced 
from a reflection group of a root system in $\mR^{n+1}\leq\mR^{n+1,1}$,
the ambient vector space of signature $(n+1,1)$. The conformal sphere $S^n$ is realized by 
projectivization of the cone $\fn$ of null length vectors in $\mR^{n+1,1}$
and is the flat version of the curved Fefferman-Graham 
ambient metric construction.    

 In the article we focus on the construction of  
differential-difference Dunkl modification of conformal Laplacian type invariants 
(also termed conformal Dunkl-Laplace operators)
associated to the structure 
$(S^n,[g_0],G)$, the conformal sphere with conformal class of the round metric and its 
finite reflection subgroup.
This collection of differential-difference operators can 
be regarded as a deformation by multiplicity function on $G$ of 
conformally invariant 
powers of the Laplace operator (corresponding to the trivial multiplicity function) 
and yields $G$-equivariant intertwining operators acting on
principal series representations of $O(n + 1, 1, \mR)$ induced from 
characters of the Levi factor of the conformal (maximal) parabolic
subgroup $P\leq O(n + 1, 1, \mR)$.   
 
In summary, the present article should be understood as a special instance 
of general research program to construct and study properties of 
differential-difference invariants
for the couple given by a manifold with geometric structure and a finite subgroup of the group 
of its automorphisms, combining techniques from the differential geometry, algebra 
and representation theory of both Lie and finite reflection groups.   


\section{Homogenous model of flat conformal geometry}

Ambient metric construction in conformal geometry associates to an $n$-dimensional conformal manifold 
$(M, [g])$ of signature $(p,q)$ a (up to certain order) pseudo-Riemannian space $(\tilde{ M},\tilde{g})$ 
of two dimensions higher, see e.g., \cite{gjms}. 

In the article we shall need just the flat version of the ambient metric 
construction, see e.g., \cite{hh}. 
The sphere $S^n$ is conformally flat manifold $(S^n, [g_0])$ 
of signature $(n, 0)$, realized as the projectivization of the cone of 
null-vectors $\fn$ in $\mR^{n+1,1}$ with the flat metric
$$
\tilde{g}(X^0,X^1,\dots ,X^n,X^\infty)=dX^0\otimes dX^\infty + dX^\infty\otimes dX^0+\sum_{i=1}^ndX^i\otimes dX^i  
$$
of signature $(n + 1, 1)$ in the coordinates $X^0,X^1,\dots ,X^n,X^\infty$.

The ambient space $\tilde{M}$ is 
a small neighborhood of $\fn$ in $\mR^{n+1,1}$ and the ambient metric is the pull-back 
of the flat metric $\tilde{g}(X^0,X^1,\dots ,X^n,X^\infty)$ from $\mR^{n+1,1}$ to
this neighborhood.  

A conformal density on $S^n$ is a section of homogeneous line bundle ${\cal L}\to S^n$, induced 
from a character of the reductive Levi factor of the conformal parabolic subalgebra 
$P$ of $O(n+1,1,\mR)$. The parabolic subalgebra $P\leq O(n+1,1,\mR)$ is the stabilizer 
of a real line in $\fn$. 

A section of ${\cal L}$ can be identified with
$\mR_+$-homogeneous function on the null-cone $\fn$ in the ambient
space. Let us remark that the construction of conformally invariant differential operators
of Laplace type acting between densities on $S^n$ proceeds by application 
of $\tilde\triangle^j$, $j\in \mN$, to an extension of a density off $\fn$ and the 
resulting section of different homogeneity restricts back to $\fn$. 

Notice that for the explicit computation of invariants in a coordinate chart and 
in the ambient metric, we need the following preparatory result not found explicitly in
\cite{gjms}, whose proof is an elementary computation.
\begin{lemma}\label{trans}
The two coordinate systems
on $\mR^{n+1,1}$: $(t,x_1,\dots ,x_n,\rho)$, 
the coordinate system adapted to the null-cone $\fn$
and $(X^0,X^1,\dots ,X^n,X^\infty)$, the euclidean coordinate system, are related by 
\begin{eqnarray}
& & t=X^0,\, x_i=\frac{X^i}{X^0},\, \rho =\frac{X^\infty}{X^0}+\frac{1}{2{X^0}^2}||X^i||^2,
\nonumber \\
& & X^0=t,\, X^i=tx_i,\, X^\infty =t(\rho -\frac{1}{2}||x_i||^2).
\end{eqnarray}
The coordinate vector fields transform as 
\begin{eqnarray}
& & \frac{\partial}{\partial X^0}=
\frac{\partial}{\partial t}-\frac{1}{t}\sum_{i=1}^nx_i\frac{\partial}{\partial x_i}
-\frac{1}{t}(\rho +\frac{1}{2}||x_i||^2)\frac{\partial}{\partial\rho},
\nonumber \\
& & \frac{\partial}{\partial X^i}=\frac{1}{t}\frac{\partial}{\partial x_i}
+ \frac{x_i}{t}\frac{\partial}{\partial\rho},\, i=1,\dots ,n
\nonumber \\
& & \frac{\partial}{\partial X^\infty}=\frac{1}{t}\frac{\partial}{\partial\rho}.
\label{coordinatechange}
\end{eqnarray}
\end{lemma}




\section{Finite reflection groups in conformal geometry 
on $S^n$ and ambient metric construction}

In what follows we restrict to the vector space $\mR^{n+1,1}$ with 
non-degenerate symmetric bilinear form $<,>$ of signature $(n+1,1)$.
If not stated otherwise, the coordinates of a vectors in $\mR^{n+1,1}$ 
are considered with respect to the canonical orthonormal basis 
$\{e_0, e_1,\dots ,e_n, e_\infty\}$.
Recall that the couple 
$(\mR^{n+1,1},<,>)$ plays the role of the flat (or homogeneous) case of the
Fefferman-Graham ambient space, \cite{gjms}. We shall adopt for our purposes
a few basic definitions
related to finite reflection groups, see e.g. \cite{bg}, \cite{ros}, \cite{eting}, \cite{opdam}, 
from euclidean signature to the signature $(n+1,1)$
and introduce the notion 
of an induced reflecting subsphere $S^{n-1}_\alpha$.   

Let us denote by $\fn\subset\mR^{n+1,1}$ the cone of null-vectors, i.e.
$\fn:=\{X\in\mR^{n+1,1}|<X,X>=0\}$. For $\alpha\in \mR^{n+1,1}\setminus\fn$ let 
$R_\alpha$ be the reflection along the hyperplane orthogonal to $\alpha$, i.e.
\begin{eqnarray}
R_\alpha X = X -2\frac{<\alpha ,X>}{<\alpha ,\alpha>}\alpha ,\, X\in\mR^{n+1,1}.
\end{eqnarray}
\begin{definition}
A finite subset $R\subset (\mR^{n+1,1}\setminus\fn )$ is called root system provided
\begin{enumerate}
\item
$R\cap\mR\alpha =\{\alpha ,-\alpha\}$,
\item
$R_\alpha(R)=R$ 
\end{enumerate}
for all $\alpha\in R$.
\end{definition} 
Let $G\leq O(\mR^{n+1,1},<,>)\simeq O(n+1,1,\mR)$ be the finite reflection group 
(generated by $\{R_\alpha|\alpha\in R\}$) associated to the root system $R$. 
We use the notation $[R]$ for the space of $G$-orbits of the root system $R$. 
An element of the vector space of all functions $k: [R]\to\mC$ is called the 
multiplicity function, and its value at $\alpha\in [R]$ is denoted $k(\alpha)$.

Note that $G$ acts linearly on $\mR^{n+1,1}$ and preserves 
$\fn$, i.e. 
\begin{eqnarray}
<X,X>=0\, \Longrightarrow\, <R_\alpha X,R_\alpha X>=0, \,\, \forall \alpha\in R.
\end{eqnarray} 
In particular, $G$ preserves the lines in $\fn$ and so induces a map on $S^n$,
the projectivization of $\fn$. The canonical (ambient) metric on $\mR^{n+1,1}$ written in 
the coordinates 
$X^0,X^1,\dots ,X^n,X^\infty$ with respect to canonical orthonormal basis, 
$$
\tilde{g}(X^0,X^1,\dots ,X^n,X^\infty)=dX^0\otimes dX^\infty + dX^\infty\otimes dX^0+\sum_{i=1}^ndX^i\otimes dX^i,  
$$
can be conveniently rewritten in the coordinates
\begin{eqnarray}
\tilde{X}^0:=\frac{1}{\sqrt{2}}(X^0+X^\infty), \, \tilde{X}^\infty:=\frac{1}{\sqrt{2}}(X^0-X^\infty),\,
\tilde{X}^i=X^i\, (i=1,\dots ,n)
\end{eqnarray}
as
$$
\tilde{g}(\tilde{X}^0,\tilde{X}^1,\dots ,\tilde{X}^n,\tilde{X}^\infty)=d\tilde{X}^0\otimes d\tilde{X}^0+
\sum_{i=1}^nd\tilde{X}^i\otimes d\tilde{X}^i - d\tilde{X}^\infty\otimes d\tilde{X}^\infty .
$$
A finite (reflection) group $G$ is contained in the maximal compact subgroup
$K=O(n,\mR)\times O(1,\mR)$ of $O(n+1,1,\mR)$, where $O(1,\mR)\simeq\mZ_2$. The 
elements in $G\cap O(n,\mR)$ are generated by a root system in the hyperplane 
$\tilde{X}^\infty=0$ (i.e. $X^0-X^\infty =0$), while the non-trivial element in
$G\cap O(1,\mR)$ is generated by reflection along the root in the coordinate axis 
$\tilde{X}^\infty$. If there was a root for $O(1,\mR)$ lying in the axis 
$\tilde{X}^\infty$, it would
have a negative norm and its reflecting hyperplane $H_\alpha$ would not intersect 
the positive null cone, 
i.e. it would not induce any fixed point on the projectivisation of the null-cone. 
In conclusion, we restrict the root system for $G$ to lie in the euclidean (with respect
to the induced metric) subspace $\tilde{X}^\infty=0$.
\begin{example}
Let us consider the irreducible root system on $\mR^{n+1}\leq\mR^{n+1,1}$ of type $B_{n+1}$,
\begin{eqnarray}
B_{n+1}:=\{\pm e_i\pm e_j)|\, 0\leq j<i\leq n\}\cup \{\pm e_i|\,  0\leq i\leq n\}.
\end{eqnarray}
The inclusion $\mR^n\leq\mR^{n+1}\leq\mR^{n+1,1}$ ($\mR^n\subset S^n$) 
implies that the set of roots $B_{n+1}$ 
splits on two subsets, $B_{n+1}=B_n\cup S$: 
\begin{enumerate}
\item
The subset $B_n$ is the set of vectors in $\mR^{n+1,1}$ 
with $\alpha_0=\alpha_{n+1}=0$,
\begin{eqnarray}
B_{n}=\{(0,\pm e_i\pm e_j,0)|1< j<i\leq n\}\cup \{(0,\pm e_i ,0)| 1< i\leq n\},
\end{eqnarray}
corresponding to standard reflections in the Levi subgroup $O(n,\mR)\leq O(n,1,\mR)$. 
\item
The subset $S$,
\begin{eqnarray}
S=\{(\pm 1,\pm e_i\pm e_j,\pm 1)|1<i\leq n\}\cup \{(1,0,\dots ,0,1)\},
\end{eqnarray}
are the roots characterized by $\alpha_0=\alpha_{n+1}\not=0$. The elements in 
$O(n+1,\mR)\setminus O(n,\mR)$ are responsible for non-trivial rational  
factor in the conformal Dunkl-Lapace operator.
\end{enumerate}
\end{example}  
So let $\alpha\in R$ be a root such that $\alpha=(\alpha_0,\alpha_1,\dots ,\alpha_n,\alpha_0)$
in canonical orthonormal basis of $\mR^{n+1,1}$. Then 
$$
<\alpha ,\alpha>=2\alpha_0^2+\sum_{i=1}^n\alpha_i^2,
$$
i.e. the fixed length $<\alpha ,\alpha>$ of $\alpha$ allows to express 
\begin{eqnarray}
\alpha_0=\sqrt{\frac{1}{2}(<\alpha ,\alpha>-\sum_{i=1}^n\alpha_i^2)}.
\end{eqnarray}
The following result is a straightforward computation.
\begin{lemma}
The reflecting hyperplane $H_\alpha\leq\mR^{n+1,1}$ associated to the root $\alpha\in R$ is  
given by linear span of $(n+1)$-tuple of (linearly independent) vectors
\begin{eqnarray}
& & (-\alpha_1,\alpha_0,0,\dots, 0),\, \nonumber \\
& & (-\alpha_2,0,\alpha_0,0,\dots, 0),\, \nonumber \\
& & \dots\, , \nonumber \\
& & (-\alpha_n,0,\dots 0,\alpha_0, 0), \nonumber \\
& & (-\alpha_0-\frac{1}{\alpha_0}\sum_{i=1}^n\alpha_i^2,\alpha_1,\dots ,\alpha_n,\alpha_0),
\end{eqnarray}
orthogonal to $\alpha$. 
\end{lemma}
\begin{definition}
The reflecting subsphere $S^{n-1}_\alpha$ corresponding to the (positive length) 
root $\alpha$ is the 
projectivization of the intersection of $H_\alpha$ and $\fn$. 
\end{definition}
In fact, in the local chart 
$U:=\{(1,x_1,\dots ,x_n,-\frac{1}{2}\sum_{i=1}^nx_i^2)|(x_1,\dots ,x_n)\in\mR^n\subset S^n\}$
of $S^n$, we have
\begin{eqnarray}
& & S^{n-1}_\alpha\cap U=\{(1,x_1,\dots ,x_n,-\frac{1}{2}\sum_{i=1}^nx_i^2)\in\fn | 
\nonumber \\
& & <(1,x_1,\dots ,x_n,-\frac{1}{2}\sum_{i=1}^nx_i^2),(\alpha_0,\alpha_1,\dots ,\alpha_n,\alpha_0)>\, =0\}.
\end{eqnarray}
The explicit form of the quadric $S^{n-1}_\alpha$ in the local chart $U$ is 
\begin{eqnarray}
\sqrt{\frac{1}{2}(<\alpha ,\alpha>-\sum_{i=1}^n\alpha_i^2)}(1-\frac{1}{2}\sum_{i=1}^nx_i^2)
+\sum_{i=1}^n\alpha_ix_i=0.
\end{eqnarray}
On the sphere $S^n$, the reflection of $(1,x_1,\dots ,x_n,-\frac{1}{2}\sum_{i=1}^nx_i^2)\in U$ along $S^{n-1}_\alpha$ is
a point on $S^n$ given by
\begin{eqnarray}\label{reflection}
X &\to & R_\alpha X = X-\frac{2<\alpha ,X>}{<\alpha ,\alpha>}\alpha , 
 \\
(1,x_1,\dots ,x_n,-\frac{1}{2}\sum_{i=1}^nx_i^2) &\mapsto &
\nonumber \\
& & \Bigl(1-\frac{2\alpha_0}{<\alpha ,\alpha>}(\alpha_0(1-\frac{1}{2}\sum_{i=1}^nx_i^2)+\sum_{i=1}^n\alpha_ix_i ),
\nonumber \\
& & x_1-\frac{2\alpha_1}{<\alpha ,\alpha>}(\alpha_0(1-\frac{1}{2}\sum_{i=1}^nx_i^2)+\sum_{i=1}^n\alpha_ix_i ),
\nonumber \\
& & \dots ,
\nonumber \\
& & x_n-\frac{2\alpha_n}{<\alpha ,\alpha>}(\alpha_0(1-\frac{1}{2}\sum_{i=1}^nx_i^2)+\sum_{i=1}^n\alpha_ix_i ),
\nonumber \\ \nonumber
& & -\frac{1}{2}\sum_{i=1}^nx_i^2-\frac{2\alpha_0}{<\alpha ,\alpha>}
(\alpha_0(1-\frac{1}{2}\sum_{i=1}^nx_i^2)+\sum_{i=1}^n\alpha_ix_i )\Bigr).
\end{eqnarray}
Let us summarize the previous considerations in the following Lemma.
\begin{lemma}\label{lemmarefl}
In the local chart $U\subset S^n$, $X\in U$, the rational map of the underlying reflection $R_\alpha$ 
associated to the root $\alpha=(\alpha_0,\alpha_1,\dots ,\alpha_n,\alpha_0)$ is
\begin{eqnarray}\label{ractrans}
& & X\to R_\alpha X,
 \\ \nonumber
& & x_i\mapsto
\frac{x_i-\frac{2\alpha_i}{<\alpha ,\alpha>}(\alpha_0(1-\frac{1}{2}\sum_{j=1}^nx_j^2)+\sum_{j=1}^n\alpha_jx_j )}
{(1-\frac{2\alpha_0}{<\alpha ,\alpha>}(\alpha_0(1-\frac{1}{2}\sum_{j=1}^nx_j^2)+\sum_{j=1}^n\alpha_jx_j )},\,\,
i=1,\dots, n.
\end{eqnarray}
\end{lemma}
Notice that Equation (\ref{reflection}) applies to all vectors $X\in\mR^{n+1,1}$, but 
the rational transformation realized by an element of the reflection group does not preserve 
the local chart $U$. So either the expression in the denominator of the rational transformation
is non-zero and the image of a given point in $U$ rests in $U$, or it is zero and the point 
maps to the point at infinity $S^n\setminus U$.
Observe that it is just the component $\alpha_0$ of the root $\alpha=(\alpha_0,\alpha_1,\dots ,\alpha_n,\alpha_0)$
responsible for the rational factor of the map $R_\alpha$.


We shall now slightly change our focus and introduce the main algebraic tool, encapsulating 
the interplay between finite reflection group $G$ and conformal geometry of $(S^n,[g_0])$. Our definition 
is basically a signature modification of certain structure in euclidean (i.e., positive definitive) signature, 
see e.g. \cite{eting},
to vector spaces of any signature, e.g. $(\mR^{n+1,1},<,>)$.
 
A basic device in the harmonic analysis for finite reflection groups $G$ on $\mR^{n+1,1}$ are 
the Dunkl differential-difference operators,
\begin{eqnarray}
& & T_\xi(k)(f)(X) := \partial_\xi f(X) +
\sum_{\alpha\in R_+}k(\alpha)
\frac{<\alpha ,\xi>}{<\alpha ,X>}
(1-R_\alpha)f(X),
\nonumber \\
& & 
\xi\in\mR^{n+1,1},\, f\in C^\infty (\mR^{n+1,1}).
\end{eqnarray}
Recall that $k$ is the multiplicity function on the root system $R$.
They form mutually commuting family of $G$-equivariant operators acting on $C^\infty(\mR^{n+1,1})$. 
Consequently, expanded in the canonical orthonormal basis $\{e_0,e_1,\dots ,e_n,e_\infty\}$ 
the Dunkl-Laplace operator
\begin{eqnarray}
\tilde\triangle_k :=\sum_{i=1}^{n+2} T_{e_i}(k)T^{e_i}(k)=\sum_{i=1}^{n+2} <T_{e_i}(k),T_{e_i}(k)>
\end{eqnarray}
fulfills $\tilde\triangle_k\circ g=g\circ \tilde\triangle_k$ for all $g\in G$.
Let us consider three homogeneous differential-difference
$G$-invariant operators in $End(C^\infty(\mR^{n+1,1}))$, written in 
signature independent notation as 
\begin{eqnarray}
& & E :=-\frac{1}{4}<X,X>,\, \nonumber \\
& & F := \tilde\triangle_k,\, \nonumber \\
& & H :=\frac{n+2}{2}+\gamma_k +<X,\partial>,
\end{eqnarray}
where $<X,\partial>$ is the Euler homogeneity operator on $\mR^{n+1,1}$ and 
$\tilde\triangle_k$ the Dunkl-Laplace operator on $\mR^{n+1,1}$:
\begin{eqnarray}\label{dunlap}
(\tilde\triangle_k f)(X)=
(\tilde\triangle f)(X) + 2\sum_{\alpha\in R_+}
k(\alpha )
(\frac
{<\tilde\nabla f (X),\alpha >}
{<\alpha ,X >}
-\frac
{f (X)-f (R_\alpha X)}
{<\alpha , X >^2 }
).
\end{eqnarray}
Here we used the notation $\tilde\triangle$ for the Laplace operator on $\mR^{n+1,1}$ and
\begin{eqnarray}
\gamma_k :=\sum_{\alpha\in R_+}k(\alpha ) .
\end{eqnarray}
The proof of the next Lemma is just a signature modification of the proof known in the 
euclidean case.  
\begin{lemma}
The three operators $E,F,G$ fulfill the $sl(2)$ commutation relations:
\begin{eqnarray}
& & [E, F] = H\, ,\nonumber \\ 
& & [E, H] = 2E\, , \nonumber \\
& & [F, H] = 2F\, .
\end{eqnarray}
\end{lemma}
The space of smooth $w$-homogeneous functions $C^\infty_w({\fam2 C})$ on the null-cone ${\fam2 C}$
is isomorphic to the space of $w$-densities $C^\infty(S^n,\kl_w)$ on $S^n$ via
\begin{eqnarray}
& & \iota: C^\infty(S^n,\kl_w)|_{\mR^n}\to C^\infty_w({\fam2 C})
\nonumber \\
& & \iota(f)(X^0,X^1,\dots ,X^n,X^\infty)=(X^0)^wf(\frac{X^1}{X^0},\dots ,\frac{X^n}{X^0}).
\end{eqnarray}
The ambient Dunkl-Laplacian $\tilde\triangle_k$ and its powers 
act on smooth $w$-densities $C^\infty(S^n,\kl_w)$ on $S^n$ as 
follows. A $w$-density on $S^n$ corresponds to a $w$-homogeneous function on 
$\fn$, which is then arbitrarily extended out of $\fn$ to its neighborhood 
in $\mR^{n+1,1}$. After the action of $\tilde\triangle_k^j$ we restrict 
back to $\fn$. For each power of $\tilde\triangle_k$ 
there exists a specific (termed critical) weight $w$ for 
which the previous procedure is independent of the extension chosen.
\begin{theorem}
\label{existence}
Let $j\in\mN$ and set $w:=-\frac{n}{2}+j-\gamma_k$.
Moreover, let $\tilde{f}$ be an extension to $\mR^{n+1,1}$ of a smooth 
$w$-density $f\in C^\infty(S^n,\kl_w)$. Then the composition
\begin{eqnarray}
f\longrightarrow\tilde{f}
\stackrel{\tilde\triangle^j_k}{\longrightarrow}
\tilde\triangle^j_k(\tilde{f})
\longrightarrow
\tilde\triangle^j_k(\tilde{f})|_\fn
\end{eqnarray}
depends on $f$ only and not on its extension $\tilde{f}$, and consequently
induces a nontrivial $G$-invariant differential-difference operator
\begin{eqnarray}
\tilde\triangle^j_k|_\fn : C^\infty(S^n,\kl_w)\longrightarrow C^\infty(S^n,\kl_{w-2j})
\end{eqnarray}
of order $2j$ and symbol $\triangle^j$.
\end{theorem}
{\bf Proof:}
Let us consider an extension $\tilde{f}$ of a $w$-homogeneous function $f$
on $\fn$. Any other extension $\tilde{f}_1$ differs from $\tilde{f}$ by $Eg$
for some $g$ of homogeneity $w-2$. We have 
$$
F\tilde{f}_1=F\tilde{f} +FEg=F\tilde{f}+EFg-(\frac{n+2}{2}+\gamma_k+w-2)g,
$$
and so for $w=-\frac{n}{2}+1-\gamma_k$ we have $F\tilde{f}_1|_\fn=F\tilde{f}|_\fn$,
exactly as claimed. The case of an integral power of $F$ is based on the 
commutator $[E,F^j]$ and is analogous to $j=1$. The proof is 
complete. 
$\blacksquare$

The operators constructed in Theorem \ref{existence} are expected to be related to 
the obstructions to harmonic extension of a $t$-homogeneous function on $\fn$ 
to a solution of the Dunkl-Laplace equation $\tilde\triangle_k \tilde{f} = 0$ on $\mR^{n+1,1}$, but
we shall not pursue this line of considerations.

It would be also extremely useful to construct the Dunkl type first order differential-difference 
operators on the sphere $S^n$, whose sum of squares yields the conformal Dunkl-Laplace 
differential-difference operators. 

\section{The explicit form of second order conformal Dunkl-Laplace operator}
In the previous section we proved the existence of Dunkl version 
of conformally invariant Laplace operator and its powers. 
An explicit form of these operators is another matter and amounts
to the computation for a chosen finite reflection group $G$ 
and in a chosen representative metric.
In the present section, we 
compute explicit form of the conformal Dunkl-Laplace operator in a local 
coordinate chart $U$ of $S^n$. 

A smooth $w$-density $f=f(x_1,\dots ,x_n)$
on $S^n$ can be regarded as a $w$-homogeneous function on $\fn$. The 
construction described in Theorem \ref{existence} is independent on the extension
of $t^wf(x_1,\dots ,x_n)$ off the null-cone $\fn$ and so we choose the constant, 
i.e. $\rho$-independent, extension, denoted by $\tilde{f}$,
\begin{eqnarray}
 \tilde{f}(t,x_1,\dots ,x_n,\rho)=\tilde{f}(t,x_1,\dots ,x_n)=t^wf(x_1,\dots ,x_n)
\end{eqnarray} 
for the critical value $w=-\frac{n}{2}-\gamma_k +1$. We determine the explicit form  
of the operator on $S^n$ induced by descending $\tilde{\triangle}_k$ (recall equation (\ref{dunlap})):
\begin{enumerate}
\item
$$
\tilde{\triangle}\tilde{f}(X)|_{t=1}=(2\frac{\partial}{\partial X^0}\frac{\partial}{\partial X^\infty}
+\sum_{i=1}^n\frac{\partial^2}{\partial {X^i}^2})\tilde{f}(X)|_{t=1}=
\sum_{i=1}^n\frac{\partial^2}{{\partial x_i}^2}f(x_1,\dots ,x_n),
$$ 
i.e. $\tilde{\triangle}$ descends in the local chart $U$ to 
the Laplace operator $\triangle=\sum_{i=1}^n\frac{\partial^2}{{\partial x_i}^2}$
acting on $f$.
\item
It follows from Lemma \ref{trans} 
$$
\tilde\nabla=((-\frac{n}{2}-\gamma_k+1-\sum_{j=1}^nx_j\frac{\partial}{\partial x_j}),
\frac{\partial}{\partial x_1},\dots ,\frac{\partial}{\partial x_n},0),
$$
so we get for $\alpha_0=\sqrt{\frac{1}{2}(<\alpha ,\alpha>-\sum_{i=1}^n\alpha_i^2)}$
\begin{eqnarray}
<\alpha ,\tilde\nabla>=\alpha_0(-\frac{n}{2}-\gamma_k+1-\sum_{j=1}^nx_j\frac{\partial}{\partial x_j})
+\sum_{i=1}^n\alpha_i\frac{\partial}{\partial x_i}.
\end{eqnarray}
Moreover,
$$
<\alpha ,X>=\sqrt{\frac{1}{2}(<\alpha ,\alpha>-\sum_{i=1}^n\alpha_i^2)}(1-\frac{1}{2}\sum_{i=1}^nx_i^2)
+\sum_{i=1}^n\alpha_ix_i,
$$
and so $\frac{<\alpha, \tilde\nabla >}{<\alpha ,X >}$ descends to the operator
\begin{eqnarray}
\frac{<\alpha ,\tilde\nabla>}{<\alpha ,X >}=
\frac{\alpha_0(-\frac{n}{2}-\gamma_k+1-\sum_{j=1}^nx_j\frac{\partial}{\partial x_j})
+\sum_{i=1}^n\alpha_i\frac{\partial}{\partial x_i}}
{\alpha_0(1-\frac{1}{2}\sum_{i=1}^nx_i^2)
+\sum_{i=1}^n\alpha_ix_i}
\end{eqnarray}
acting on $f$.
\item
It follows from the formula (\ref{ractrans}) and comments beyond Lemma \ref{lemmarefl}
\begin{eqnarray}
\tilde{f}(X)|_{t=1}&=&f(x_1,\dots ,x_n),
\nonumber \\
\tilde{f}(R_\alpha X)|_{t=1}&=&\tilde{f}\biggl(
X^0-\frac{2\alpha_0}{<\alpha ,\alpha>}(\alpha_0(X^0-\frac{1}{2X^0}\sum_{i=1}^n{X^i}^2)+\sum_{i=1}^n\alpha_iX^i),
\nonumber \\
& & X^1-\frac{2\alpha_1}{<\alpha ,\alpha>}(\alpha_0(X^0-\frac{1}{2X^0}\sum_{i=1}^n{X^i}^2)+\sum_{i=1}^n\alpha_iX^i),
\nonumber \\
& & \dots ,
\nonumber \\
& & X^n-\frac{2\alpha_n}{<\alpha ,\alpha>}(\alpha_0(X^0-\frac{1}{2X^0}\sum_{i=1}^n{X^i}^2)+\sum_{i=1}^n\alpha_iX^i),
\nonumber \\
& & -\frac{1}{2X^0}\sum_{i=1}^n{X^i}^2-\frac{2\alpha_0}{<\alpha ,\alpha>}
(\alpha_0(X^0-\frac{1}{2X^0}\sum_{i=1}^n{X^i}^2)+\sum_{i=1}^n\alpha_iX^i )\biggr)|_{t=1}
\nonumber \\
&=& \biggl(1-\frac{2\alpha_0}{<\alpha ,\alpha>}
(\alpha_0(1-\frac{1}{2}\sum_{i=1}^n{x_i}^2)+\sum_{i=1}^n\alpha_ix_i)\biggr)^{-\frac{n}{2}-\gamma_k +1}
\nonumber \\
& & f\biggl(
\frac{x_1-\frac{2\alpha_1}{<\alpha ,\alpha>}(\alpha_0(1-\frac{1}{2}\sum_{i=1}^n{x_i}^2)+\sum_{i=1}^n\alpha_ix_i)}
{1-\frac{2\alpha_0}{<\alpha ,\alpha>}(\alpha_0(1-\frac{1}{2}\sum_{i=1}^n{x_i}^2)+\sum_{i=1}^n\alpha_ix_i)},
\nonumber \\
& & \dots ,
\nonumber \\
& & 
\frac{x_n-\frac{2\alpha_n}{<\alpha ,\alpha>}(\alpha_0(1-\frac{1}{2}\sum_{i=1}^n{x_i}^2)+\sum_{i=1}^n\alpha_ix_i)}
{1-\frac{2\alpha_0}{<\alpha ,\alpha>}(\alpha_0(1-\frac{1}{2}\sum_{i=1}^n{x_i}^2)+\sum_{i=1}^n\alpha_ix_i)}
\biggr)
\end{eqnarray}
and so
\begin{eqnarray}
\frac{(1-R_\alpha)\tilde f(X)}{<\alpha ,X >^2}|_{t=1}&=&\frac{\tilde f(X)-\tilde f(R_\alpha X)}{<\alpha ,X >^2}|_{t=1}=
\nonumber \\
& & \frac{1}{(\alpha_0(1-\frac{1}{2}\sum_{i=1}^nx_i^2)+\sum_{i=1}^n\alpha_ix_i)^2}\cdot \biggl(f(x_1,\dots ,x_n)-
\nonumber \\
& & \Bigl(1-\frac{2\alpha_0}{<\alpha ,\alpha>}(\alpha_0(1-\frac{1}{2}\sum_{i=1}^n{x_i}^2)+\sum_{i=1}^n\alpha_ix_i)\Bigr)^{-\frac{n}{2}-\gamma_k +1}\cdot
\nonumber \\
& & f\Bigl(\frac{x_1-\frac{2\alpha_1}{<\alpha ,\alpha>}(\alpha_0(1-\frac{1}{2}\sum_{i=1}^n{x_i}^2)+\sum_{i=1}^n\alpha_ix_i)}
{1-\frac{2\alpha_0}{<\alpha ,\alpha>}(\alpha_0(1-\frac{1}{2}\sum_{i=1}^n{x_i}^2)+\sum_{i=1}^n\alpha_ix_i)},
\nonumber \\
& & \dots ,
\nonumber \\
& & 
\frac{x_n-\frac{2\alpha_n}{<\alpha ,\alpha>}(\alpha_0(1-\frac{1}{2}\sum_{i=1}^n{x_i}^2)+\sum_{i=1}^n\alpha_ix_i)}
{1-\frac{2\alpha_0}{<\alpha ,\alpha>}(\alpha_0(1-\frac{1}{2}\sum_{i=1}^n{x_i}^2)+\sum_{i=1}^n\alpha_ix_i)}\Bigr)\biggr).
\nonumber \\
\end{eqnarray}
\end{enumerate}
\begin{theorem}
Let $G\leq O(n+1,\mR)\leq O(n+1,1,\mR)$ be a finite reflection group associated to the root system
$R=\{\alpha |\alpha=(\alpha_0,\alpha_1,\dots ,\alpha_n,\alpha_0)\}\subset\mR^{n+1,1}$. Then 
the conformal Dunkl-Laplace operator is in the local chart $U\simeq\mR^n\subset S^n$ given by
\begin{eqnarray}
& & \triangle_kf(x_1,\dots ,x_n)=\sum_{i=1}^n\frac{\partial^2}{{\partial x_i}^2}f(x_1,\dots ,x_n)+ 2\sum_{\alpha\in R_+}k(\alpha)
\nonumber \\
& & 
\Biggl(
\frac{1}{(\alpha_0(1-\frac{1}{2}\sum_{i=1}^nx_i^2)+\sum_{i=1}^n\alpha_ix_i)}\cdot
\nonumber \\
& & \Bigl(\alpha_0(-\frac{n}{2}-\gamma_k+1-\sum_{j=1}^nx_j\frac{\partial}{\partial x_j})
+\sum_{i=1}^n\alpha_i\frac{\partial}{\partial x_i}\Bigr)f(x_1,\dots ,x_n)
\nonumber \\
& & -\frac{1}{(\alpha_0(1-\frac{1}{2}\sum_{i=1}^nx_i^2)+\sum_{i=1}^n\alpha_ix_i)^2}\cdot 
\nonumber \\
& & \biggl(f(x_1,\dots ,x_n)-
\nonumber \\
& & \Bigl(1-\frac{2\alpha_0}{<\alpha ,\alpha>}(\alpha_0(1-\frac{1}{2}\sum_{i=1}^n{x_i}^2)+\sum_{i=1}^n\alpha_ix_i)\Bigr)^{-\frac{n}{2}-\gamma_k +1}\cdot
\nonumber \\
& & f\Bigl(\frac{x_1-\frac{2\alpha_1}{<\alpha ,\alpha>}(\alpha_0(1-\frac{1}{2}\sum_{i=1}^n{x_i}^2)+\sum_{i=1}^n\alpha_ix_i)}
{1-\frac{2\alpha_0}{<\alpha ,\alpha>}(\alpha_0(1-\frac{1}{2}\sum_{i=1}^n{x_i}^2)+\sum_{i=1}^n\alpha_ix_i)},
\nonumber \\
& & \dots ,
\nonumber \\
& & 
\frac{x_n-\frac{2\alpha_n}{<\alpha ,\alpha>}(\alpha_0(1-\frac{1}{2}\sum_{i=1}^n{x_i}^2)+\sum_{i=1}^n\alpha_ix_i)}
{1-\frac{2\alpha_0}{<\alpha ,\alpha>}(\alpha_0(1-\frac{1}{2}\sum_{i=1}^n{x_i}^2)+\sum_{i=1}^n\alpha_ix_i)}\Bigr)\biggr)\Biggr),
\end{eqnarray}
when acting on $f\in C^\infty (S^n,{\lL}_w)$ for the critical weight $w$. $\blacksquare$
\end{theorem}
The explicit formulas for conformal powers of the second order Dunkl-Laplace operator are
notationally less trivial due to the iterative composition of rational terms. 
Because of complexity of the number of descending terms we do not attempt to write down 
extensive formulas explicitly.


\section{Open questions}

 There are several interesting aspects of the construction, leading 
to its generalization in various directions. First of all, the present
construction works in any signature. Also, one does not need to restrict
to the principal series representations of $O(n + 1, 1, \mR)$ induced from 
character, but 
inducing from a finite dimensional representation (e.g. the spinor
representation) leads to the canonical $G$-equivariant version of a 
given conformally invariant operator (e.g. the conformally invariant 
powers of the Dirac operator.) Another direction for possible generalization 
stems from the existence of a curved version of the Fefferman-Graham ambient
metric construction of the manifold with conformal structure $(M,[g])$,
see the review \cite{gjms} and the references therein. 
In particular, the ambient metric is an Einstein 
metric of constant scalar curvature, and so
finite reflection subgroups of the (compact subgroup of) the  
isometry group $Isom(\tilde{M},\tilde{[g]})$ of the Fefferman-Graham metric 
of signature $(n+1,1)$ together with its
action on the Fefferman-Graham ambient space can be geometrically analyzed 
along the lines of e.g., \cite{akm}. Similar considerations might be 
applied to finite reflection groups realized in the group of automorphisms
of another geometric structure, e.g. the complex reflection groups leading to the 
Dunkl version of so called CR-Laplace operator and its powers, see \cite{gg}.
Another source of examples involves quaternion reflection groups and related
parabolic geometries of quaternion type. 

 On general abstract grounds, the problems discussed in the article are rather special 
instances of the classification scheme for $G$-invariant differential-difference 
operators on the sphere with standard round metric and its conformal structure 
$(S^n, [g_0])$, related to 
the branching problems for a (infinite dimensional) representations 
of a Lie group and finite subgroups $G\leq O(n + 1, 1, \mR)$. To our best knowledge, 
there is rather poor understanding of such problems. Moreover, in addition to the 
deformation of operators of elliptic type one should attempt to construct
the Dunkl version of overdetermined (twistor type) operators and more 
generally Dunkl deformation of the geometric construction of Bernstein-Gelfand-Gelfand 
sequences, \cite{css}, or the invariant calculus of tractors for a
wide class of parabolic geometries, \cite{cg}.
Another generalization of the present results includes the variation on the theme of finite reflection
groups. There is no need for such restrictive category, i.e. one can consider for 
a fixed background geometry a wider category of subgroups like the class of 
finite groups or even discrete groups.
 
 The author expects that compatible structure consisting of a finite reflection group
and the underlying conformal structure might be fruitful for the development 
of a conformal (or CR, quaternion, etc.) version of invariant theory for 
(complex, quaternion) reflection groups accompanied by algebraic structures
of integrable systems, Cherednik and Hecke algebras, etc., in a systematic way.   

The Fefferman-Graham ambient space $(\tilde{M},\tilde{g})$ associated to a general 
conformal manifold $(M,[g])$ is a pseudo-Riemannian manifold equipped with an Einstein 
metric, i.e. a particularly nice class of constant scalar curvature manifolds on 
which the action of a finite refection group can be studied.    

The following definition is a natural generalization of a finite refection group acting 
on unitary vector space.  
\begin{definition}
A reflection $s$ on a Riemannian manifold $(N,g)$ is an isometry, $s\in Isom(N)$, 
such that for some fixed point $x$ of $s$ the tangent map $T_xs$ is a reflection 
in the Hilbert space $(T_xN,g_x)$.  
\end{definition}
For a discrete subgroups of $Isom(N)$ generated by reflections, the classical 
geometrical concepts like its Dirichlet domain, Weyl chamber and for Coxeter 
manifolds a Riemannian chamber defined as a manifold with corners such that 
its walls are totally geodesic submanifolds and neighboring walls satisfy the Coxeter 
property, are studied, see e.g. \cite{akm} and the references therein. 
The techniques of Fefferman-Graham ambient metric construction are sufficiently
flexible allowing in particular instances of finite 
reflection groups to construct Dunkl deformations of GJMS-operators or more 
generally, differential-difference invariants of $(M,[g],G)$ on 
manifolds with conformal structure $(M,[g])$. However, at the moment it is
not clear what are the discrete subgroups of $Isom(\tilde{M},\tilde{g})$
of most interest worth to be studied in details.


\flushleft{{\em Acknowledgment}: 
The author gratefully acknowledges the support of the grant GA CR P201/12/G028. 


\vspace{0.5cm}
\flushleft{Petr Somberg\\
Mathematical Institute of Charles University,\\
Sokolovsk\'a 83, 186 75 Praha\\ Czech Republic,\\
somberg@karlin.mff.cuni.cz}

\end{document}